\documentclass[a4letter, 12pt,reqno]{amsart} 
\newif\ifdraft
\draftfalse

\ifdraft\usepackage[notref,notcite]{showkeys}\fi
\usepackage{graphicx,subfigure,suffix}

\usepackage{amssymb,amsmath,amsfonts,amssymb,mathtools}
\usepackage{latexsym}
\usepackage{cite}
\usepackage[pdfborder={0 0 .1}]{hyperref}

\makeatletter
\@mparswitchfalse
\makeatother

\begin{document}

 \title [\resizebox{4.5 in}{!}{Fast growth of the vorticity gradient in symmetric smooth domains for 2D incompressible ideal flow}]{ Fast growth of the vorticity gradient in symmetric smooth domains for 2D incompressible ideal flow}
\begin{abstract}
We construct an initial data for the two-dimensional Euler equation in a bounded smooth symmetric domain such that the gradient of vorticity in $L^{\infty}$ grows as a double exponential in time for all time. Our construction is based on the recent result by Kiselev and \v{S}ver\'{a}k \cite{kiselev2014small}.
\end{abstract}

  \author{Xiaoqian Xu}
  \address{\hskip-\parindent
  Xiaoqian Xu\\
  Department of Mathematics\\
  University of Wisconsin-Madison\\
  Madison, WI 53706, USA}
  \email{xxu@math.wisc.edu}

\maketitle
\section{Introduction}
The two dimensional Euler equation in vorticity form for an incompressible fluid is given by
\begin{equation}
\partial_t\omega+(u\cdot \nabla) \omega=0,\quad \omega(x,0)=\omega_0(x).
\end{equation}
Here $\omega=\mbox{curl}\  u$ is the vorticity of the flow, and the velocity $u$ can be determined from $\omega$ by the Biot-Savart law. Here we consider the flow in a smooth bounded domain $\Omega$, and we assume $u$ satisfies a no flow boundary condition, namely $u\cdot n=0$ on $\partial \Omega$, here $n$ is the outer normal vector of $\partial \Omega$. This implies that
$$
u=\nabla^{\perp}\int_{\Omega}G_{\Omega}(x,y)\omega(y) dy.
$$
Where $G_{\Omega}(x,y)$ is the Green's function of the Dirichlet problem in $\Omega$ and $\nabla^{\perp}=(\partial_{x_2},-\partial_{x_1})$ is the perpendicular gradient.\\

Global regularity of solutions to 2D Euler equation is well-known, it was first proved in \cite{Wolibner1933} and \cite{Holder1933}, and see also \cite{kato1967classical}, \cite{marchioro1994mathematical}, \cite{chemin1998perfect} for more work. By the standard estimates (see \cite{yudovich1963flow}), to obtain the higher regularity of the solution corresponding to smooth initial data, we only need to bound the $L^{\infty}$ norm of the gradient of vorticity. The best known upper bound for this quantity has double exponential growth in time. This result is well-known and has first appeared in \cite{yudovich1963flow}. \\

The sharpness of the double exponential growth bound has remained open for a long time. In \cite{yudovich2000loss}, Yudovich provided an example showing unbounded growth of the vorticity gradient at the boundary of the domain. Nadirashvili \cite{nadirashvili1991wandering} constructed an example on the annulus where the gradient of vorticity grows linearly in time. 
In recent years, Denisov did a series of work on this problem. In \cite{denisov2009infinite}, he constructed an example on the torus with superlinear growth in time. In \cite{denisov2015double}, he proved that double exponential growth can be possible for any fixed but finite time. In \cite{denisov2015sharp}, he constructed a patch solution to smoothly forced 2D Euler equation where the distance of two patches decreases double exponentially in time. We refer to \cite{kiselev2014small} and \cite{zlatovs2015exponential} for more information on these questions.\\

In 2014, Kiselev and \v{S}ver\'{a}k \cite{kiselev2014small} constructed an example of initial data in the disk where  such double exponential growth was observed. This means the double exponential growth in time is the best upper bound we can get for the gradient of vorticity, at least for solutions to 2D Euler equation in the disk. How generic is such growth 
is an interesting question. If this growth happens in many situations, it would highlight the challenges involved into numerical simulation of the solutions. In this paper, we 
make the first step towards better understanding of this question, by generalizing the construction of Kiselev and Sverak to the case of an arbitrary sufficiently regular domain 
with a symmetry axis. The main new ingredient of our proof is a more general version of the key lemma in \cite{kiselev2014small}. The lemma captures hyperbolic structure of the velocity 
field near the point on the boundary where the fast growth happens. The analysis in \cite{kiselev2014small} is specific to the disk, and the Green's function estimates used 
to prove the lemma do not translate easily to the more general case. We develop a more flexible natural construction that allows to carry the estimates through. \\

In this paper,  we will prove the following double exponential growth bound for more general domains instead of disk:\\
\textbf{Theorem 1.} Let $\Omega$ be a $C^3$ bounded open domain in $\mathbb{R}^2$, tangent at the origin to the $x_1$-axis and symmetric about $x_2$-axis. Consider the 2D Euler equation on $\Omega$. There exists a smooth initial data $\omega_0$ with $\|\nabla\omega_0\|_{L^{\infty}}>\|\omega_0\|_{L^{\infty}}$ such that the corresponding solution $\omega(x,t)$ satisfies
\begin{equation}\label{mainin}
\dfrac{\|\nabla\omega(\cdot,t)\|_{L^{\infty}}}{\|\omega_0\|_{L^{\infty}}}\geq \left(\dfrac{\|\nabla\omega_0\|_{L^{\infty}}}{\|\omega_0\|_{L^{\infty}}}\right)^{c(\Omega)\exp(c(\Omega)\|\omega_0\|_{L^{\infty}}t)}
\end{equation}
for some $c(\Omega)>0$ that only depends on $\Omega$ and for all $t\geq 0$.\\

Based on the ideas of \cite{kiselev2014small}, we will focus on the appropriate representation for the Biot-Savart law for the fluid velocity $u$. 
The key will be obtaining the representation which will show that $u_1\sim x_1\log(x_1)$ on appropriate (time dependent) length scales. 
After that, based on the construction of \cite{kiselev2014small}, we will get the desired growth for the gradient of vorticity in domain $\Omega$.

\section{The key lemma}
To construct a flow with fast growth in the gradient, we need a technical lemma for the expansion of velocity field. \\

\noindent
We use the same notation as in \cite{kiselev2014small}, which means $\omega$ is the vorticity field in $\Omega$, odd in the first variable. Note that this property is 
conserved by Euler evolution in a symmetric domain. Let $u$ be the corresponding velocity field. $D^+=\{x\in \Omega: x_1>0\}, \tilde{x} = (-x_1,x_2)$ and $Q(x_1,x_2)$ is a region that is the intersection of $D^+$ and the quadrant $\{y:x_1\leq y_1<\infty$, $x_2\leq y_2<\infty\}$.\\

\noindent
\textbf{Key Lemma.} Suppose $\Omega$ is a $C^3$ bounded open domain in $\mathbb{R}^2$, symmetric about $x_2$-axis, and tangent to $x_1$-axis at the origin. Take any $\gamma$, $\frac{\pi}{2}>\gamma>0$. Denote $D_1^{\gamma}$ the intersection of $D^+$ with a sector $\frac{\pi}{2}-\gamma\geq \phi\geq -\frac{\pi}{2}$, where $\phi$ is the usual angular variable. Then there exists $\delta>0$ such that for all $x\in D_1^{\gamma}$ such that $|x|\leq \delta$ we have
\begin{equation}\label{u1}
u_1(x,t)=-\frac{4}{\pi}x_1\int_{Q(x_1,x_2)}\dfrac{y_1y_2}{|y|^4}\omega(y,t)dy+x_1B_1(x,t),
\end{equation}
where $|B_1(x,t)|\leq C(\gamma,\Omega)\|\omega_0\|_{L^{\infty}}$.\\
Similarly, if we denote $D_2^{\gamma}$ the intersection of $D^+$ with a sector $\frac{\pi}{2}\geq \phi\geq \gamma$, then for all $x\in D^{\gamma}_2$ such that $|x|\leq \delta$, we have
\begin{equation}\label{mainth}
u_2(x,t) = \frac{4}{\pi}x_2\int_{Q(x_1,x_2)}\dfrac{y_1y_2}{|y|^4}\omega(y,t)dy+x_2B_2(x,t).
\end{equation}
where $|B_2(x,t)|\leq C(\gamma,\Omega)\|\omega_0\|_{L^{\infty}}$.\\

In the argument below, $\delta$ and constant $C$ may change from line to line, but since we change it only for finite many times, at the end we choose the smallest $\delta$ and the biggest $C$ among all of them instead.\\

\noindent
We want to define $y^*$ to be the mirror image point of $y$, about the boundary. Namely, we want $y^*$ to be written as $y^*=2e(y)-y$, where $e(y)$ is a point on $\partial\Omega$ so that $y-e(y)$ is orthogonal to the tangent line at $e(y)$. More intuitively, $e(y)$ is the closest point to $y$ on $\partial \Omega$. However, it is not clear that this $e(y)$ is well-defined. Given $y$, it may be possible to find more than one $e(y)$ on the whole $\partial\Omega$. However, we will show $y^*(y)$ is locally well-defined close to the origin by \textbf{Lemma 1} below.\\

\noindent
Since $\partial \Omega$ is tangent to the $x_1$ axis at the origin, we can choose the parameterization near the origin of $\partial\Omega$ such that we have $\partial \Omega=(s,f(s))$ for some function $f\in C^3$ and sufficiently small $s.$\\
\noindent
\textbf{Lemma 1.} There exists $r=r(\Omega)>0$ only depending on $\Omega$, so that for any $y\in B_{r}(0)\cap \Omega$, there is a unique $s$ in $(-2r,2r)$ such that the following equation holds:
\begin{equation}\label{min}
(s-y_1)+f'(s)(f(s)-y_2)=0.
\end{equation}
Moreover, this $s$, as a function of $y$, is $C^2$.\\
\textbf{Remark.} If we call $e(y)=(s(y),f(s(y)))$, then $(\ref{min})$ means $y-e(y)$ is orthogonal to the tangent line of $\partial \Omega$ at $e(y)$.
Note that if $e(y)$ is one of the points on $\partial \Omega$ closest to $y$ then \eqref{min} holds. \\

\begin{proof}[Proof of lemma 1.] We call the left hand side of $(\ref{min})$ the function $F(s,y)$. We take derivative of $F(s,y)$ in $s$ and get
\begin{equation}\label{implicit}
1+f'(s)^2+f''(s)(f(s)-y_2).
\end{equation}
So $D_sF(0,(0,0))=1$ and $F(0,(0,0))=0$. Thus by implicit function theorem we know the solution to $(\ref{min})$ exists and is unique in a neighborhood of $\{0\}\times\{(0,0)\}$. By choosing $r$ small enough, $s(y)$ is uniquely defined. Moreover, as in the very beginning we choose $\partial\Omega$ to be $C^3$, which means $f$ is $C^3$. Hence, $F(s,y)$ is a $C^2$ function as it only contains function $f'$ and $f$, which means $s(y)$ is $C^2$ in $y$ by the implicit function theorem.
\end{proof}
\noindent
By \textbf{Lemma 1}, we can give the following definition:\\

\noindent
\textbf{Definition 1.} Given $r$ small enough, for any $|y|\leq r$, we define $e(y)$ to be the only point in $B_{2r}(0)\cap \partial\Omega$ so that $e(y)-y$ is orthogonal to the tangent line of $\partial\Omega$ at $e(y)$. And we define $y^*(y) = 2e(y)-y$.  We denote $y^*(y)=(y^*_1(y),y^*_2(y))=(y^*_1,y_2^*)$.\\

\noindent
Then, we have the following lemma.\\
\noindent
\textbf{Lemma 2.} Take $r$ small enough as in \textbf{Lemma 1}. If $|s_0|\leq 2r$ and $y\in B_{r}(0)\cap \Omega$, we have
\begin{equation}\label{expansion0}
\begin{split}
&y^*_1-s_0=\frac{1-f'(s_0)^2}{1+f'(s_0)^2}(y_1-s_0)+\frac{2f'(s_0)}{1+f'(s_0)^2}(y_2-f(s_0))+O(|y-(s_0,f(s_0))|^2),\\
&y^*_2-f(s_0)=
\frac{2f'(s_0)}{1+f'(s_0)^2}(y_1-s_0)+\frac{f'(s_0)^2-1}{1+f'(s_0)^2}(y_2-f(s_0))+O(|y-(s_0,f(s_0))|^2).
\end{split}
\end{equation}
Here the constant in capital O notation only depends on the domain $\Omega$. In addition the map $y^*$ is invertible in $y$ in $B_{r}(0)\cap \Omega$.\\

\begin{proof}[Proof of lemma 2.] This is an elementary calculation by Taylor's expansion formula. Like in the proof of \textbf{Lemma 1}, we call the left hand side of $\eqref{min}$ $F(s,y)$. We take the partial derivative of $F(s,y)$ with $s=s(y)$ in $y_1$ and denote $\partial_1=\partial_{y_1}$, $\partial_2=\partial_{y_2}$. We get
$$
\partial_1s-1+f'(s)^2\partial_1s+f''(s)\partial_1 s(f(s)-y_2)=0.
$$
Plug in $y_1=s_0$ ,$y_2=f(s_0)$, noticing that we have $s(s_0,f(s_0))=s_0$, we get
$$
\partial_1s|_{y=(s_0,f(s_0))}=\frac{1}{1+f'(s_0)^2}.
$$
Which means by the definition of $y^*$ we get
$$
\partial_1y_1^*|_{y=(s_0,f(s_0))}=2\partial_1s|_{y=(s_0,f(s_0))}-1=\frac{1-f'(s_0)^2}{1+f'(s_0)^2}
$$
And similarly by taking the partial derivative in $y_2$ in $F(s,y)$ we get
$$
 \partial_2 s|_{y=(s_0,f(s_0))}=\frac{f'(s_0)}{1+f'(s_0)^2}.
$$
So
$$
\partial_2y_1^*|_{y=(s_0,f(s_0))}=2\partial_2s|_{y=(s_0,f(s_0))}=\frac{2f'(s_0)}{1+f'(s_0)^2}
$$
By chain rule, we get
$$
\partial_1f(s)|_{y=(s_0,f(s_0))}= \frac{f'(s_0)}{1+f'(s_0)^2},\quad \partial_2f(s)|_{y=(s_0,f(s_0))}=\frac{f'(s_0)^2}{1+f'(s_0)^2}.
$$
Which means we can get
$$
\partial_1y_2^*|_{y=(s_0,f(s_0))}=\frac{2f'(s_0)}{1+f'(s_0)^2},\quad\partial_2y_2^*|_{y=(s_0,f(s_0))}=\frac{f'(s_0)^2-1}{1+f'(s_0)^2}.
$$
Thus by Taylor's expansion we get $(\ref{expansion0})$. And we can also see the invertibility of $y^*$ for $r$ small enough, this is simply by the inverse function theorem because $|\det(\nabla y^*)|_{y\in\partial\Omega}|=1$.
\end{proof}

To understand $y^*$ better, we need another lemma. The following lemma shows that although $\partial\Omega$ could be crazy, the intuition that $y^*$ must be outside of $\Omega$ is always true.\\

\noindent
\textbf{Lemma 3.} There exists $r$ so that for any $y\in B_r(0)\cap \Omega$, $y^*(y)\notin \Omega$.

\noindent
\begin{proof}[Proof of lemma 3.] First, since $y$ is close to the origin, the slope of inner normal line at $e(y)$ of $\partial\Omega$ is close to $+\infty$. Recall by our definition the inner normal line has the same direction as $y-e(y)$, which means the second component of $e(y)$ is less than $y_2$. By the definition of $y^*$, $y^*_2<y_2$.\\
\noindent
Now we argue by contradiction. Suppose for every $r_0$ we can find $y$ so that $y^*(y)$ is inside $\Omega$. By the expansion of $y^*$ near zero, we know that $|y^*|\approx |y|$. Here the notation $"\approx"$ means there are constants $C_1,C_2$ only depending on $\Omega$, such that $C_1|y|\leq |y^*|\leq C_2|y|$. So, if $r_0$ is small enough, say, less than $\frac{r}{C_2}$, where $r$ is the same $r$ in \textbf{Lemma 1}, then $y^*(y)$ is also in the domain of map $y^*$ by \textbf{Lemma 1}. By definition, $e(y)-y^*(y)$ is orthogonal to the tangent line at $e(y)$. However, by \textbf{Lemma 1}, such a boundary point $e(y)$ is unique, which means $e(y)=e(y^*(y))$. So we know $y^*(y^*(y))=y$. But then $y_2=y^*(y^*(y))_2<y^*_2<y_2$, which is a contradiction.
\end{proof}

We will use $y^*$ as a sort of conjugate point for $y$ in the context of the Dirichlet reflection principle for the representation of the Green's function.
We note that for the case of the disk in \cite{kiselev2014small}, by the well known explicit formula for the Green's function,
the natural choice of $y^*$ is given by circular inversion of $y$. For more general $\Omega$ the choice of $y^*$ is less obvious. We will see that our definition of $y^*$ will work well for the estimates that we have in mind.\\

\noindent
Without loss of generality, we assume $\Omega\subset [-2,2]\times [-2,2]$. Then we have the following proposition.\\

\noindent
\textbf{Proposition 1.} Suppose $\Omega$ is a $C^3$ bounded open domain in $\mathbb{R}^2$, symmetric about $x_2$-axis, and tangent to $x_1$-axis at the origin. Then there exists $r=r(\Omega)>0$ so that for $x,y\in B_{r}(0)$, the Green function of $\Omega$ can be written as:
\begin{equation}\label{expansion}
G_{\Omega}(x,y)=\frac{1}{2\pi}(\log|x-y|-\log|x-y^*|)+B(x,y).
\end{equation}
Here $B(x,y)$ satisfies for any $\omega\in L^{\infty}(\Omega)$, $\int_{B_{r}(0)\cap\Omega}B(x,y)\omega(y)dy \in C^{2,\alpha}(B_{r}(0)\cap \Omega)$, for any $0<\alpha<1$. More precisely, we have 
$$
\|\partial_{x_i}\partial_{x_j}\int_{B_{r}(0)\cap\Omega}B(x,y)\omega(y)dy\|_{L^{\infty}(B_{r}(0)\cap \Omega)}\leq C(\Omega)\|\omega\|_{L^{\infty}}\quad i,j=1,2
$$
\\

\noindent
To prove this proposition, we need a technical lemma.\\

\noindent
\textbf{Lemma 4.} Let $x=(s,f(s))$. Let $K(z_1,z_2)$ be a integral kernel such that it is $C^1$ on the set $\{z_1\in \Omega, z_2\in \Omega:z_1\neq z_2\}$. Suppose we have $K(s,f(s),y)=O(\frac{1}{|x-y|})$ and $D_sK(s,f(s),y)=O(\frac{1}{|x-y|^2})$. Then $\int_{\Omega} K(s,f(s),y)\omega(y)dy $ has modulus of continuity $s\log(s)$, with the constant equal to $C(\Omega)||\omega||_{L^{\infty}}$. Here $C(\Omega)$ is a constant that only depends on $\Omega$.\\

\begin{proof}[Proof of lemma 4.] Without loss of generality, let $s_1<s_2$. Suppose $|s_1-s_2|=\zeta$. Let $(\frac{(s_1+s_2)}{2},f(\frac{(s_1+s_2)}{2}))$ be $Z$. By the smoothness of $f$, there is a constant $C_1=C_1(\Omega)$, so that for any $s$ between $s_1$ and $s_2$, $(s,f(s))\in B_{C_1 \zeta}(Z)$. Then, for any $\tau>C_1 \zeta$, we have
\begin{equation}\label{ball}
\begin{split}
&\int_{\Omega} (K(s_1,f(s_1),y)-K(s_2,f(s_2),y))\omega dy \\
&=\int_{\Omega\cap B_{\tau}(Z)} (K(s_1,f(s_1),y)-K(s_2,f(s_2),y))\omega(y) dy+\\
                                           & \int_{\Omega\cap B^c_{\tau}(Z)} (K(s_1,f(s_1),y)-K(s_2,f(s_2),y))\omega(y) dy\\
                                           &\leq C ||\omega||_{L^{\infty}}\int_{\Omega\cap B_{\tau}(Z)}(\frac{1}{|y-(s_1,f(s_1))|}+\frac{1}{|y-(s_2,f(s_2))|})dy+\\
&  C ||\omega||_{L^{\infty}}\int_{\Omega\cap B^c_{\tau}(Z)} \int_{s_1}^{s_2} |D_sK(t,f(t),y)|dt dy\\
\end{split}
\end{equation}
For $y\in B_{\tau}(Z)$, by the smoothness of $f$,
$$
|y-(s_i,f(s_i))|\leq |y-Z|+|Z-(s_i,f(s_i))|\leq \tau+C_1 \zeta,
$$
for $i=1,2$. In addition, for any $s_1\leq s\leq s_2$ and $y\in B^c_{\tau}(Z)$, we have
$$
|y-(s,f(s))|\geq |y-Z|-|Z-(s,f(s))|\geq \tau-C_1 \zeta.
$$
Hence, the right hand side of $(\ref{ball})$ is no more than
\begin{equation}
\begin{split}
&C ||\omega||_{L^{\infty}}\int_{(\Omega-Z)\cap B_{\tau+C_1 \zeta}(O)}\frac{1}{|y|}dy+C ||\omega||_{L^{\infty}}\int_{(\Omega-Z)\cap B^c_{\tau-C_1 \zeta}(O)} \frac{|s_1-s_2|}{|y|^2} dy \leq \\
& C ||\omega||_{L^{\infty}}(\int_0^{\tau+C_1 \zeta}\frac{1}{r}\cdot rdr+\int_{\tau-C_1 \zeta}^2\frac{\zeta}{r^2}\cdot rdr)\\
                                           &=C||\omega||_{L^{\infty}}(\tau+\zeta\log(\tau-C_1 \zeta)+C_1 \zeta).
\end{split}
\end{equation}
Here $\Omega-Z$ means the translation of $\Omega$ by $Z$. So if we choose $\tau= 4C_1\zeta$, we get the desired modulus of continuity.
\end{proof}
\noindent
\textbf{Remark.} In this lemma, it's easy to see that if $K(s,f(s),y)$ is not differentiable but $|K(s_1,f(s_1),y)-K(s_2,f(s_2),y)|=|s_1-s_2|^{\alpha}O(\frac{1}{|(s_1,f(s_1))-y|^{2}}+\frac{1}{|(s_2,f(s_2)-y)|^2})$, we can still get the similar result. More precisely, $\int_{\Omega} K(s,f(s),y)\omega(y)dy$ has modulus of continuity $x^{\alpha}\log(x)$. This can be used to extend the results of this paper to less regular domains with $C^{2,\alpha}$ boundary. \\

\noindent
Now we prove \textbf{Proposition 1}.\\
\noindent
\begin{proof}[Proof of Proposition 1]

The idea is to use the properties of elliptic equations. 

First, remember 
$$
B(x,y)=G_{\Omega}(x,y)-\frac{1}{2\pi}(\log|x-y|-\log|x-y^*|),
$$
where $G_{\Omega}(x,y)$ is the Green function of domain $\Omega$, $y^*$ is a function of $y$ defined by \textbf{Definition 1}. As a well-known result, Green function is a smooth function for $x\neq y$. Here we would like to show the subtraction of $\frac{1}{2\pi}\log|x-y|$, which is the Green function of $\mathbb{R}^2$, can eliminate the singularity of $G_{\Omega}(x,y)$ with the help of the term $\log|x-y^*|$. 

As a result of \textbf{Lemma 3}, we know that for all $y$, $y^*\notin \Omega$. Therefore for any fixed  $y\in\Omega$, $\log |x-y^*|$ is smooth and harmonic in $x$. This means $B(x,y)$ is harmonic as $x$ varies in $\Omega$, and satisfies the boundary condition $B(x,y)|_{x\in\partial\Omega}=\frac{1}{2\pi}\log(\frac{|x-y^*|}{|x-y|})$. Hence $\int_{\Omega\cap B_{r}(0)}B(x,y)\omega(y)dy$ is also harmonic and satisfies $$\int_{\Omega\cap B_{r}(0)}B(x,y)\omega(y)dy|_{x\in\partial\Omega\cap B_{r}(0)}=\int_{\Omega\cap B_{r}(0)}\frac{1}{2\pi}\log(\frac{|x-y^*|}{|x-y|})\omega(y)dy.$$ 
Since the boundary of the domain $\Omega$ is $C^3,$ by the well-known results on elliptic regularity (see, e.g., Lemma $6.18$ of \cite{gilbarg2001elliptic}), we know that in order to show that $\int_{\Omega\cap B_{r}(0)}B(x,y)\omega(y)dy$ is $C^{2,\alpha}$ near the origin, we only need to show that this harmonic function is $C^{2,\alpha}$ on the boundary near the origin. Recall the notation $x=(s,f(s))$. We only need
\begin{equation}\label{boundary}
\iota(s)=\int_{\Omega\cap B_{r}(0)}\frac{1}{2\pi}\log(\frac{|(s,f(s))-y^*|}{|(s,f(s))-y|})\omega(y)dy
\end{equation}
to be $C^{2,\alpha}$ in $s$ for $s$ small. Here remember $y^*$ is only a function in $y$, so $\iota(s)$ is a well defined function in $s$. The proof for regularity of $\iota$ is simply by calculation. Here we will only use the expansion of $x-y^*$ and the corresponding cancellation of $y-y^*$. As $\iota(s)$ can be seen as the integral of a difference of the same function in different points $y$ and $y^*$, we will essentially need to calculate the finite differences of some complicated functions.  

First we find the second derivative in $s$ of $\log(\frac{|(s,f(s))-y^*|}{|(s,f(s))-y|})$. We call it $K(s, f(s),y)$. More precisely, we have
\begin{equation}
\begin{split}
-K(s,\mu,y)&=\frac{1+f''(s)(\mu-y_2)+f'(s)^2}{|x-y|^2}-2\frac{((s-y_1)+f'(s)(\mu-y_2))^2}{|x-y|^4}\\
       &-\frac{1+f''(s)(\mu-y^*_2)+f'(s)^2}{|x-y^*|^2}+2\frac{((s-y^*_1)+f'(s)(\mu-y^*_2))^2}{|x-y^*|^4}.
\end{split}
\end{equation}
Then, observe that by \textbf{Lemma 2} and simple computation, for $y$ close to $(s,f(s))$ we have $f(s)-y_2^*=f(s)-y_2+O(|f(s)-x_2|)+O(|s-x_1|)+O(|x-y|^2)=f(s)-y_2+O(|x-y|)$, which means $y_2-y_2^*=O(|x-y|)$ for $y$ close to $x$, and $|x-y^*|^2=|x-y|^2+O(|x-y|^3)$, for $x$, $y$ close to the origin and $y$ close to $x$. So $K(s,f(s),y)$ can be written as
\begin{equation}\label{KXY}
\begin{split}
 &\frac{f''(s)(f(s)-y_2+O(|x-y|))+O(|x-y|))}{|x-y|^2+O(|x-y|^3)}-\\
 &2\frac{A(s,y)(s-y_1+s-y^*_1+f'(s)(f(s)-y_2+f(s)-y^*_2))}{|x-y|^4+O(|x-y|^5)}+\\
 &2\frac{((s-y_1)+f'(s)(f(s)-y_2))^2O(|x-y|)}{|x-y|^4+O(|x-y|^5)}.
\end{split}
\end{equation}
Where $A(s,y)=(s-y_1)-(s-y^*_1)+f'(s)((f(s)-y_2)-(f(s)-y^*_2)))$. By \textbf{Lemma 2},
\begin{equation}
\begin{split}
A(s,y)&=(s-y_1)+\frac{1-f'(s)^2}{1+f'(s)^2}(y_1-s)+\frac{2f'(s)}{1+f'(s)^2}(y_2-f(s))+\\
&f'(s)\left(f(s)-y_2+\frac{2f'(s)}{1+f'(s)^2}(y_1-s)+\frac{f'(s)^2-1}{1+f'(s)^2}(y_2-f(s))\right)+\\
&O(|x-y|^2)\\
&=\frac{-2f'(s)^2}{1+f'(s)^2}(y_1-s)+\frac{2f'(s)}{1+f'(s)^2}(y_2-f(s))+\\
&f'(s)\left(\frac{2f'(s)}{1+f'(s)^2}(y_1-s)+\frac{-2}{1+f'(s)^2}(y_2-f(s))\right)+O(|x-y|^2)\\
&=O(|x-y|^2).
\end{split}
\end{equation}
Again by \textbf{Lemma 2} we have $s-y^*_1=O(|x-y|)$, $f(s)-y^*_2=O(|x-y|)$. And also we have $s-y_1=O(|x-y|)$ and $f(s)-y_2=O(|x-y|)$. Plug in all of these into $(\ref{KXY})$ we get $K(s,f(s),y)=O(\frac{1}{|x-y|})$.\\
\noindent
Then, we take the derivative of $K(s,f(s),y)$ in terms of $s$ again. We get
\begin{equation}
\begin{split}
D_sK(s,f(s),y)&=\frac{f'''(s)(f(s)-y_2)+3f''(s)f'(s)}{|x-y|^2}-\\
&\frac{(1+f''(s)(f(s)-y_2)+f'(s)^2)((s-y_1)+f'(s)(f(s)-y_2))}{|x-y|^4}\\
&+4(\frac{1+f''(s)(f(s)-y_2)+f'(s)^2}{|x-y|^2}-2\frac{((s-y_1)+f'(s)(f(s)-y_2))^2}{|x-y|^4})\\
&\cdot \frac{((s-y_1)+f'(s)(f(s)-y_2))}{|x-y|^2}\\
&-\frac{f'''(s)(f(s)-y^*_2)+3f''(s)f'(s)}{|x-y^*|^2}\\
&+\frac{(1+f''(s)(f(s)-y^*_2)+f'(s)^2)((s-y^*_1)+f'(s)(f(s)-y^*_2))}{|x-y^*|^4}\\
&-4(\frac{1+f''(s)(f(s)-y^*_2)+f'(s)^2}{|x-y^*|^2}-2\frac{((s-y^*_1)+f'(s)(f(s)-y^*_2))^2}{|x-y^*|^4})\\
&\cdot \frac{((s-y^*_1)+f'(s)(f(s)-y^*_2))}{|x-y^*|^2}
\end{split}
\end{equation}
\begin{equation}
\begin{split}
&= O(\frac{1}{|x-y|^2})\\
&-5(1+f'(s)^2)(\frac{(s-y_1+f'(s)(f(s)-y_2))}{|x-y|^4}-\frac{(s-y^*_1+f'(s)(f(s)-y^*_2))}{|x-y^*|^4})\\
&-8(\frac{((s-y_1)+f'(s)(f(s)-y_2))^3}{|x-y|^6}-\frac{((s-y^*_1)+f'(s)(f(s)-y^*_2))^3}{|x-y^*|^6})\\
&=O(\frac{1}{|x-y|^2})-5(1+f'(s))\frac{A(s,y)}{|x-y|^4}-8\frac{A(s,y)O(|x-y|^2)}{|x-y|^6}\\
&=O(\frac{1}{|x-y|^2}).
\end{split}
\end{equation}
To complete the proof of this proposition, we only need to use \textbf{Lemma 4}.
\end{proof}

\textbf{Remark 1.} Notice that this proposition is true for all small enough $r$. Later in the proof of the key lemma this $r$ may change from line to line, and finally we will choose the smallest $r$ which is still only depend on $\Omega$.\\

\textbf{Remark 2.} If $f$ is not in $C^3$ but in $C^{2,\beta}$ for some $0<\beta<1$, by a longer but similar computation we can find
$$
|K(s_1,f(s_1),y)-K(s_2,f(s_2),y)|=|s_1-s_2|^{\alpha}O(\frac{1}{|(s_1,f(s_1))-y|^2}+\frac{1}{|(s_2,f(s_2))-y|^2}).
$$
Which means even if we have $C^{2,\beta}$ domain, we can still get some regularity of $K(s,f(s),y)$. By the remark after \textbf{Lemma 4}, we still have $\int_{B_{r}(0)\cap\Omega}B(x,y)\omega(y)dy$ is $C^{2,\alpha}$ for any $0<\alpha <\beta$.\\

\noindent
The proposition can now be applied to prove the key lemma for the domain $\Omega$.\\

\noindent
\begin{proof}[Proof of the key lemma.]
\noindent
By \textbf{Proposition 1}, we know we can write the Green function of $\Omega$ as follows:
\begin{equation}\label{expression}
2\pi G_{\Omega}(x,y)=\eta_{B_{r}(0)}(y)(\log|x-y|-\log|x-y^*|)+C(x,y).
\end{equation}
Here $\eta_{B_{r}}(y)$ is the smooth cut-off function. $C(x,y)$ is a function so that $\int_{\Omega} C(x,y)\omega(y) dy$ is $C^{2,\alpha}(B_{\delta}(0)\cap \Omega)$, for any small $\delta\leq \frac{r}{2}$, and $\omega(y)$ is a bounded function in $\overline{\Omega}$. Here $y^*$ is the same as in \textbf{Proposition 1}. Hence, by using the Taylor's expansion and $x_2$ can be controlled by $x_1$ in the sector $D^{\gamma}_1$, with $\frac{|x|}{x_1}\leq C(\gamma)$, the first order term of $\partial_{x_2}\int_{\Omega} C(x,y)\omega(y) dy$ can be written as $x_1 J_1(x,t)+M_1(\omega)$, for $J_1(x,t)\leq C(\gamma)||\omega_0||_{L^{\infty}}$ and $M_1(\omega)=\partial_{x_2}\int_{\Omega}C(x,y)\omega(y)dy|_{x=(0,0)}$.\\

\noindent
We first prove $(\ref{u1})$. For $(\ref{mainth})$, it is similar. By the expansion of $y^*$ near the origin we know $|y^*|\approx |y|$ for any $y\in B_r(0)$. Without loss of generality, we assume $|y^*|\geq C_1 |y|$ for some $C_1\leq 1$. Fix a small $\gamma>0$, fix $x\in D_1^{\gamma}$, $|x|\leq \delta$. Here $\delta$ is a small number that we will choose later. Now we would like to choose a number which is comparable to $x_1$ while it can control both $x_1$ and $x_2$. We define $a=\frac{100}{C_1}(1+\cot(\gamma))x_1$. Let's first assume $\delta$ is small enough so that $a<\mbox{min}\{0.01,
\frac{r}{2}\}$ whenever $|x|\leq \delta$. Now the contribution to $u_1$ from integration over $B_a(0)$ does not exceed
\begin{equation}
\begin{split}
|2\pi\int_{\Omega\cap B_a(0)}\partial_{x_2}G_{\Omega}(x,y)\omega(y)dy|\leq C\|\omega_0\|_{L^{\infty}}\int_{D^+\cap B_a(0)}\left(\frac{1}{|x-y|}+1\right)dy\\
\leq C a \|\omega_0\|_{L^{\infty}}\leq C(\gamma) x_1 \|\omega_0\|_{L^{\infty}}.
\end{split}
\end{equation}
For $y\in (D^+\cap B_r(0))\setminus B_a(0)$, we have $|y|\geq 100|x|$ and $|y^*|\geq 100|x|$. By symmetry, we can write the first term in $(\ref{expression})$ as $\eta_{D^+\cap B_r(0)}$ times the following terms:
\begin{equation}\label{log}
\begin{split}
\log|x-y|-\log|x-y^*|=&\log\left(1-\frac{2xy}{|y|^2}+\frac{|x|^2}{|y|^2}\right)-\log\left(1-\frac{2xy^*}{|y^*|^2}+\frac{|x|^2}{|y^*|^2}\right)\\
&-\log\left(1-\frac{2\tilde{x}y}{|y|^2}+\frac{|x|^2}{|y|^2}\right)+\log\left(1-\frac{2\tilde{x}y^*}{|y^*|^2}+\frac{|x|^2}{|y^*|^2}\right).
\end{split}
\end{equation}
Here $\tilde{x}=(-x_1,x_2)$. For small $t$, we have
$$
\log(1+t)=t-\frac{t^2}{2}+O(t^3).
$$
Hence, $(\ref{log})$ can be written as
$$
-\frac{x_1y_1}{|y|^2}+\frac{x_1y_1^*}{|y^*|^2}-\frac{2x_1x_2y_1y_2}{|y|^4}+\frac{2x_1x_2y^*_1y^*_2}{|y^*|^4}+O(\frac{|x|^3}{|y|^3}).
$$
In the last term, we used that $|y^*|\approx |y|$, this is true by taking $s_0=0$ in the expression of $y^*$ in \textbf{Lemma 2} for $y\in B_{r}(0)$ and $r$ small. Again by the expression near $0$ of $y^*$, we have
$$
\frac{y_1^*}{|y^*|^2}=\frac{y_1+O(|y|^2)}{|y|^2+O(|y|^3)}=\frac{y_1}{|y|^2+O(|y|^3)}+b_1(y)=\frac{y_1}{|y|^2}+b_1(y).
$$
Where $b_1(y)$ is a bounded function in $y$, and the bound is a universal constant. Similarly,
$$
\frac{y^*_2}{|y^*|^2}=-\frac{y_2}{|y|^2}+b_2(y).
$$
Here again $b_2(y)$ is also bounded by a universal constant. Therefore we get that the expression $(\ref{log})$ can be written as
\begin{equation}\label{logg}
x_1b_1(y)-\frac{4x_1x_2y_1y_2}{|y|^4}+\frac{2x_1x_2y_1}{|y|^2}b_2(y)+\frac{2x_1x_2y_2}{|y|^2}b_1(y)+O(\frac{|x|^3}{|y|^3}).
\end{equation}
Then we can differentiate the above expression with respect to $x_2$, since we know the explicit functions, and by direct computation we get
$$
-\frac{4x_1y_1y_2}{|y|^4}+\frac{2x_1y_1}{|y|^2}b_2(y)+\frac{2x_1y_2}{|y|^2}b_1(y)+O(\frac{|x|^2}{|y|^3}).
$$
Now
$$
\int_{(D^+\cap B_r(0))\setminus B_a(0)}\frac{|x|^2}{|y|^3}dy\leq C|x|^2\int_a^1\frac{1}{s^2}ds\leq C\frac{|x|^2}{a}\leq C(\gamma) x_1.
$$
Also,
$$
\int_{(D^+\cap B_r(0))\setminus B_a(0)}\frac{y_i}{|y|^2}dy\leq C,
$$
for $i=1,2$.
Therefore, the last three terms of $(\ref{logg})$ only give regular contributions to $u_1$. Now we only need to show that adjusting the region $B_r(0)\setminus B_a(0)$ to $Q(x_1,x_2)$ will not change too much for the expression, namely,
$$
\int_{(D^+\cap B_r(0))\setminus B_a(0)} \frac{y_1y_2}{|y|^4}\omega(y)dy=C(\Omega)b_3(x)\|\omega_0\|_{L^{\infty}}+\int_{Q(x_1,x_2)} \frac{y_1y_2}{|y|^4}\omega(y)dy.
$$
Here again $b_3(x)$ is a bounded function whose bound is a universal constant. Indeed,
$$
\int_{D^+\setminus B_r(0)}\frac{y_1y_2}{|y|^4}\omega(y)dy\leq C(r)\|\omega_0\|_{L^{\infty}}\leq C(\Omega)\|\omega_0\|_{L^{\infty}}.
$$
And
\begin{equation}
\begin{split}
&\left|\int_{B_a\cap Q(x_1,x_2)}\frac{y_1y_2}{|y|^4}\omega(y)dy\right|\leq C\|\omega_0\|_{L^{\infty}}\int_{B_a\cap Q(x_1,x_2)}\frac{y_1|y_2|}{|y|^4}dy\\
&\leq C\|\omega_0\|_{L^{\infty}}2\int_{x_1}^{Cx_1}dy_1\int_0^{Cx_1}dy_2\frac{y_1y_2}{|y|^4}\leq C\|\omega_0\|_{L^{\infty}}.
\end{split}
\end{equation}
Finally, the set $D^+\setminus (B_a\cup Q(x_1,x_2))$ consists of two strips. The contribution of the strip along $x_2$ axis does not exceed the following quantity:
\begin{align*}
\left|\int_{D^+\setminus (B_a\cup Q(x_1,x_2))\cap \{y_1\leq x_1\}}\frac{y_1y_2}{|y|^4}\omega(y)dy\right|\leq \|\omega_0\|_{L^{\infty}}\int_{0}^{x_1}dy_1\int_{x_1}^1dy_2\frac{y_1y_2}{|y|^4}\\
\leq \|\omega_0\|_{L^{\infty}}\int_0^{Cx_1}\frac{y_1}{C^2x_1^2+y_1^2}dy_1\leq C\|\omega_0\|_{L^{\infty}}.
\end{align*}
Similarly, the integral over the strip along $x_1$ axis can be bounded by
\begin{align*}
\left|\int_{D^+\setminus (B_a\cup Q(x_1,x_2))\cap \{y_2\leq x_2\}}\frac{y_1y_2}{|y|^4}\omega(y)dy\right|\leq \left|\int_{-C(\Omega)|x_1|^2}^{|x_2|}dy_2\int_{Cx_1}^1dy_1\frac{y_1y_2}{|y|^4}\right|\|\omega_0\|_{L^{\infty}}\\
\leq C\|\omega_0\|_{L^{\infty}}\int_0^{C(\gamma)x_1}dy_2\int_{Cx_1}^1dy_1 \frac{y_1y_2}{|y|^4}\leq C(\gamma)\|\omega_0\|_{L^{\infty}}.
\end{align*}
Here the first step is due to the fact if we write $\partial \Omega=(s,f(s))$, then since $f'(0)=0$, near $0$ we have $f(x_1)=Cx_1^2$. The second inequality is true since $\delta$ is small, $|x_1|^2\leq |x|^2\leq |x|\leq C(\gamma)x_1$. This completes the estimate of the first term of $(\ref{expression})$.\\

\noindent
Finally notice that $u_1(0,0)=0$, so $M_1(\omega)$ will be canceled by the constant term of the first term. So we finish the proof of the key lemma.\\
\end{proof}
\noindent
\textbf{Remark.} By the remark after \textbf{Lemma 4} and \textbf{Proposition 1}, one can find that this key lemma is still true for $\partial \Omega$ to be $C^{2,\alpha}$, for any $\alpha>0$. Therefore one could have double exponential in time upper bound as well. On the other hand, it has been proved in \cite{kiselev2015blow} and \cite{itoh2014remark} that if the boundary $\partial\Omega$ is only Lipschitz, one may expect finite time blowup or exponential in time upper bound for $\|\nabla\omega\|_{L^{\infty}}$. It is an interesting question whether we could get any similar estimate for $C^{1,\alpha}$ domain.\\

\section{The proof of the main theorem}
Now based on the key lemma for $\Omega$, we follow the idea of the proof in \cite{kiselev2014small}, we can prove \textbf{Theorem 1}.\\
\begin{proof}[Proof of Theorem 1.]  First of all, we set $1\geq \omega_0\geq 0$ with $\|\omega_0\|_{L^{\infty}}=1$. Then we know $1\geq \omega\geq 0$ as well. If $x_2\leq 0$, observe that
\begin{equation}\label{diff}
\left|\int_{x_1}^2\int_{x_2}^{-x_2}\dfrac{y_1y_2}{|y|^4}\omega(y)dy_2dy_1\right|\leq C\int_{x_1}^2\int_{0}^{-f(x_1)}\dfrac{y_1y_2}{|y|^4}dy_1dy_2\leq C \log(1+(\frac{f(x_1)}{x_1})^2)+C\leq C.
\end{equation}
So if we take smooth $\omega_0$ equal to one everywhere in $D^+$ except on a thin strip of width $\delta$ near the axis $x_1=0$, where $0\leq \omega_0\leq 1$, we will have
$$
\int_{Q(x_1,x_2)}\dfrac{y_1y_2}{|y|^4}\omega(y)dy_1dy_2\geq C_1\int_{2\delta}^2\int_{\frac{\pi}{6}}^{\frac{\pi}{3}}\dfrac{\omega(r\cos\phi,r\sin\phi)}{r}d\phi dr-C,
$$
here in the second inequality we set $y_1=r\cos\phi,y_2=r\sin\phi$. Since $\omega<1$ in $D^+$ only in an area not exceeding $2\delta$, for $\delta$ small enough, the right hand side will be at least
\begin{equation}\label{logd}
\frac{C_1}{2}\int_{\delta}^2\int_{\frac{\pi}{6}}^{\frac{\pi}{3}} \dfrac{1}{r}dr\geq c\log(\delta^{-1}),
\end{equation}
for some $c>0.$ 

\noindent
For $0<x'_1<x''_1<1$ we denote
\begin{equation}
R(x'_1,x''_1)=\{(x_1,x_2)\in D^+, x'_1<x_1<x''_1,x_2<x_1\}.
\end{equation}
For $0<x_1<1$ we define
\begin{equation}
u^l_1(x_1,t)=\min_{(x_1,x_2)\in D^+, x_2<x_1} u_1(x_1,x_2,t)
\end{equation}
and
\begin{equation}
u^u_1(x_1,t)=\max_{(x_1,x_2)\in D^+, x_2<x_1} u_1(x_1,x_2,t).
\end{equation}
By the smoothness of $u$, it is easy to see that these functions are locally Lipschitz in $x_1$, with the Lipschitz constant bounded in finite time. As a result, we can define $a(t)$ and $b(t)$ by
\begin{equation}
\begin{split}
&\dot{a}=u_1^u(a,t),\quad a(0)=\epsilon^{10},\\
&\dot{b}=u_1^l(b,t),\quad b(0)=\epsilon.\\
\end{split}
\end{equation}
Where $0<\epsilon<\delta$ is small and to be determined later. Let $R_t=R(a(t),b(t))$. Notice by definition $R_t$ can only be guaranteed to be non-empty for small enough $t$, but we will see that in fact $R_t$ is not empty for all $t>0$.\\
We assume $\omega_0=1$ on $R_0$ and smoothly become $0$ in the $\epsilon^{10}$-neighborhood of $R_0$.  Our claim is in $R_t$ $\omega$ is always $1$ for $\delta$ small enough. \\

By the \textbf{Key Lemma} and $(\ref{diff})$, we know that $u_1$ is negative for small $\delta$. Hence both $a(t)$ and $b(t)$ are decreasing functions of time. And by $(\ref{logd})$, near the diagonal $x_1=x_2$ for $|x|<\delta$ we have
\begin{equation}\label{dia}
\frac{x_1(\log(\delta^{-1})-C)}{x_2(\log(\delta^{-1})+C)}\leq \frac{-u_1(x_1,x_2)}{u_2(x_1,x_2)}\leq \frac{x_1(\log(\delta^{-1})+C)}{x_2(\log(\delta^{-1})-C)}.
\end{equation}
This means that the vector field $u$ is directed out of the region $R_t$ on the diagonal.
In addition, by the definition of $a(t)$ and $b(t)$, the fluid trajectories starting at the points outside of $R_0$ cannot enter $R_t$ at any positive time through the vertical segments $\{(a(t),x_2)\in D^+,x_2<a(t)\}$ and $\{(b(t),x_2)\in D^+,x_2<b(t)\}$. Therefore, trajectories originating outside $R_0$ will not enter $R_t$ at any time. 
This means that $\omega=1$ in $R_t$.\\
We call $\Lambda(x_1,x_2,t)=\frac{4}{\pi}\int_{Q(x_1,x_2)}\frac{y_1y_2}{|y|^4}\omega(y)dy_1dy_2$. By the \textbf{Key Lemma} we have
$$
u^l_1(b(t),t)\geq -b(t)\Lambda(b(t),x_2(t))-Cb(t).
$$
If $x_2\leq 0$, then $x_2\geq f(x_1)\geq -Cx_1^2 $. Otherwise if $x_2>0$, $x_2\leq x_1$. By an estimate similar to (\ref{diff}) and the fact $x_2\leq b(t)$ in $R(t)$, we know
\begin{align*}
|\Lambda(b(t),x_2(t))|&\leq |\Lambda(b(t),b(t))|+\left|\frac{4}{\pi}\int_{b(t)}^2\int_{x_2}^{b(t)}\frac{y_1y_2}{|y|^4}dy_1dy_2\right|\\
&\leq |\Lambda(b(t),b(t))|+\left|\frac{4}{\pi}\int_{b(t)}^2\int_{f(b(t))}^{b(t)}\frac{y_1y_2}{|y|^4}dy_1dy_2\right|\\
&\leq |\Lambda(b(t),b(t))|+\left|\frac{4}{\pi}\int_{b(t)}^2\int_{0}^{b(t)}\frac{y_1y_2}{|y|^4}dy_1dy_2\right|\\
&\leq  |\Lambda(b(t),b(t))| +\left|\frac{2}{\pi}\int_{b(t)}^2y_1\left(\frac{1}{y^2_1}-\frac{1}{y_1^2+b(t)^2}\right)dy_1 \right|\\
&\leq  |\Lambda(b(t),b(t))| +C ,
\end{align*}
for some constant $C\geq 0$. Therefore we get
\begin{equation}\label{44}
u^l_1(b(t),t)\geq -b(t)\Lambda(b(t),b(t))-Cb(t).
\end{equation}
And by a similar estimate we also have
$$
u^u_1(a(t),t)\leq -a(t)\Lambda(a(t),0)+Ca(t).
$$
Observe that by geometry of the regions involved we also have
$$
\Lambda(a(t),0)\geq \frac{4}{\pi}\int_{R_t}\dfrac{y_1y_2}{|y|^4}dy_1dy_2+\Lambda(b(t),b(t)).
$$
Since $\omega=1$ on $R_t$ and if $\epsilon$ is sufficiently small,
$$
\int_{R_t}\dfrac{y_1y_2}{|y|^4}dy_1dy_2\geq \int_{\frac{\pi}{100}}^{\frac{\pi}{4}}\int_{\frac{a(t)}{\cos\phi}}^{\frac{b(t)}{\cos\phi}}\dfrac{\sin 2\phi}{2r}dr d\phi\geq C(-\log a(t)+\log b(t))-C.
$$
As a result,
\begin{equation}\label{45}
u^u_1(a(t),t)\leq -a(t)\left(-C(\log a(t)-\log b(t))+\Lambda(b(t),b(t))\right)+Ca(t).
\end{equation}
 Then from the estimates (\ref{44}) and (\ref{45}) we know $a(t)$ and $b(t)$ are monotone decaying in time, and by finiteness of $\|u\|_{L^{\infty}}$ these function are Lipschitz in $t$. So we have sufficient regularity to do the following calculations
$$
\dfrac{d}{dt}\log(b(t))\geq -\Lambda(b(t),b(t))-C,
$$
and
$$
\dfrac{d}{dt}\log(a(t))\leq C(\log(a(t))-\log(b(t)))-\Lambda(b(t),b(t))+C.
$$
Hence, by subtraction we have
$$
\dfrac{d}{dt}(\log a(t)-\log b(t))\leq C(\log a(t)-\log b(t))+2C.
$$
By Gronwall's inequality we get $\log a(t)\leq (9\epsilon+C)\exp(\frac{t}{C})$, and by choosing $\epsilon$ small enough, we have $a(t)\leq \epsilon^{8\exp(C t)}$. Note that the first coordinate of the characteristic originating at the point on $\partial\Omega$ near the origin with $x_1=\epsilon^{10}$, does not exceed $a(t)$ by the definition of $a(t)$. To get (\ref{mainin}), we only need to choose the initial data $\omega_0$ such that $\|\nabla\omega_0\|\lesssim \epsilon ^{-10}$. Thus, by the mean value theorem applying to $\omega$ between the origin and the point $(a(t),a(t))$, we get the desired lower bound with $\|\omega_0\|_{L^{\infty}}=1$.
\end{proof}
\noindent
\textbf{Acknowledgment.} The author acknowledges the support of NSF grants 1104415 and 1159133. I would like to thank Prof. Alexander Kiselev for the introduction into the subject, and guidance. I would also like to thank Sigurd Angenent, Mikhail Feldman and Tam Do for helpful discussions. I would like to thank the referees for providing many useful suggestions.

\end{document}